\documentclass[11pt]{article}
\newtheorem{theorem}{Theorem}
\newtheorem{lemma}[theorem]{Lemma}
\newtheorem{prop}[theorem]{Proposition}

 \newcommand{\nl}{\newline}
 \newcommand{\dist}{{\rm dist}}
 \newcommand{\N}{{\bf N}}
 
\newcommand{\R}{{\bf R}}

 \newcommand{\diam}{{\rm diam}}

\newcommand{\ia}{({\rm i})}
\newcommand{\ib}{({\rm ii})}

\newcommand{\sscr}[2]{\scriptstyle {#1} \atop\scriptstyle {#2}}

\title{Improved Rellich inequalities for the polyharmonic operator}

\author{
G. Barbatis
}

\date{}

\begin{document}

\maketitle

\begin{abstract}
We prove two improved versions of the Hardy-Rellich inequality for the polyharmonic
operator $(-\Delta)^m$ involving the distance to the boundary. The first involves an
infinite series improvement using logarithmic functions, while the second contains $L^2$ norms
and involves as a coefficient the volume of the domain. We find explicit constants for
these inequalities, and we prove their optimality in the first case.

\

\noindent {\bf AMS 2000 MSC:} 35J20 (35P99,
26D10, 47A75)
\nl {\bf Keywords:}  Polyharmonic operator, Hardy-Rellich inequality, distance to the boundary, sharp constants
\end{abstract}

\section{Introduction}\label{sec1}

Let $\Omega$ be a convex domain in $\R^n$ and let $d(x)=\dist(x,\partial\Omega)$. The classical Hardy's inequality
asserts that
\begin{equation}
\int_{\Omega}|\nabla u|^2 dx \geq\frac{1}{4}\int_{\Omega}\frac{u^2}{d^2}dx \; , \qquad u\in C^{\infty}_c(\Omega).
\label{hi}
\end{equation}
There has recently been an increased interest in so-called inproved Hardy's inequalities, where additional non-negative terms
appear in the right-hand side of (\ref{hi}). Such inequalities were first established by Maz'ya \cite{M} in the case where $\Omega$
is a half-space. 
Renewed interest in such inequalities followed the work of Brezis and Marcus
\cite{BM} where (\ref{hi}) was improved in two ways. More precisely, let
$X_1(s)=(1-\log s)^{-1}$, $s\in (0,1]$, a function that vanishes at logarithmic speed at $s=0$.
It is shown in \cite{BM} that if $\Omega$ is bounded with diameter 
$D$ then there holds
\begin{equation}
\int_{\Omega}|\nabla u|^2 dx \geq\frac{1}{4}\int_{\Omega}\frac{u^2}{d^2}dx 
+\frac{1}{4}\int_{\Omega}\frac{u^2}{d^2}X^2(d/D)dx\; , \qquad u\in C^{\infty}_c(\Omega),
\label{hi1}
\end{equation}
and also
\begin{equation}
\int_{\Omega}|\nabla u|^2 dx \geq\frac{1}{4}\int_{\Omega}\frac{u^2}{d^2}dx 
+\frac{1}{4D^2}\int_{\Omega}\frac{u^2}{d^2}dx\; , \qquad u\in C^{\infty}_c(\Omega).
\label{hi2}
\end{equation}
Inequalities (\ref{hi1}) and (\ref{hi2}) subsequently led to additional improvements and generalizations, which broadly can be termed
logarithmic and non-logarithmic respectively.

Let us define recursively $X_i(s) = X_{1}(X_{i-1}(s))$, $i\geq2$, $s\in (0,1]$.
Hence the $X_i$'s are iterated logarithmic functions that vanish at an increasingly slow rate at $s=0$ and satisfy $X_i(1)=1$.
In was proved in \cite{BFT2} that for any $p>1$ there exists $D\geq \sup_{\Omega}d(x)$ such that
\begin{equation}
\int_{\Omega}|\nabla u|^pdx \geq \Big(\frac{p-1}{p}\Big)^p\int_{\Omega}\frac{|u|^p}{d^p}dx
+  \frac{1}{2}\Big(\frac{p-1}{p}\Big)^{p-1}\sum_{i=1}^{\infty}
 \int_{\Omega}\frac{|u|^p}{d^p}X_1^2(d/D) \ldots X_i^2(d/D) dx ,
\end{equation}
for all $u\in C^{\infty}_c(\Omega)$. Each new term in this series is optimal, with respect to both the exponent two of $X_i$
and the constant $(1/2)((p-1)/p)^{p-1}$. An analogous result for the bilaplacian is obtained in \cite{BT} where it is shown that
\begin{equation}
\int_{\Omega}(\Delta u)^2dx \geq \frac{9}{16}\int_{\Omega}\frac{u^2}{d^4}dx
+  \frac{5}{8}\sum_{i=1}^{\infty}
 \int_{\Omega}\frac{u^2}{d^2}X_1^2(d/D) \ldots X_i^2(d/D) dx ,
\label{bt}
\end{equation}
which is, again, sharp.

Concerning non-logarithmic inequalities and answering a question of \cite{BM},
Hoffmann-Ostenhof et al. \cite{HHL} proved that $\diam(\Omega)^{-2}$ in (\ref{hi2})
can be replaced by $c|\Omega|^{-2/N}$, where $|\Omega|$ stands for the volume of $\Omega$; more precisely, they showed that
\begin{equation}
\int_{\Omega}|\nabla u|^2dx\geq\frac{1}{4}\int_{\Omega}\frac{u^2}{d^2}dx
 +\frac{N}{4}\Big(\frac{|\Omega|}{a_N}\Big)^{-2/N}\int_{\Omega}u^2dx,
\label{hhl}
\end{equation}
where, here and below, $a_N$ stands for the volume of the unit ball in $\R^N$. This was generlized
to $p\neq 2$ by Tidblom \cite{T} who obtained
\begin{equation}
\int_{\Omega}|\nabla u|^pdx\geq\Big(\frac{p-1}{p}\Big)^p
\int_{\Omega}\frac{|u|^p}{d^p}dx +
(p-1)\Big(\frac{p-1}{p}\Big)^p\frac{\sqrt{\pi}\Gamma(\frac{N+p}{2})}{\Gamma(\frac{p+1}{2})\Gamma(\frac{N}{2})}
 \Big(\frac{a_N}{|\Omega|} \Big)^{\frac{p}{N}}\int_{\Omega}|u|^pdx .
\label{tid}
\end{equation}
Such inequalities, where the volume of $\Omega$ appears in the right-hand side, have also been called {\em geometric}, and we follow this terminology.
In the case of geometric improvements the identification
of best constants is significantly more complex, since the problem has a global character
as opposed to local in the logarithmic case. Results in this direction where obtained in \cite{BFT3} in the linear case and when
$\Omega$ is the unit ball $B$; in particular, the best constant was identified in dimension $N=3$. The constants appearing in (\ref{hhl}) and
(\ref{tid}) are not sharp. A different type of non-logarithmic $L^p$ improvemnts, rather in the spirit of \cite{M}, is obtained in \cite{T2}.
See also \cite{FMT1,FMT2} for recent results on improved $L^p$ Hardy-Sobolev inequalities, where an $L^q$ norm, $q>p$,
is added to the right-hand side of Hardy's inequality.

The Hardy-Rellich inequalities have various applications in the study of
elliptic and parabolic PDE's.
Improved Rellich inequalities are useful if critical potentials
are additionally present and they also serve to identify such potentials. As the simplest example, one obtains
information on the existence of solution and asymptotic behavior for
the equation $u_t=\Delta+V$ (or $u_t=-\Delta^2+V$) for critical potentials $V$. We refer to \cite{D,BM,O,MMP,BT} and references therein
for more on applications.

Our aim in this article is the study of analogous problems for the polyharmonic operator $(-\Delta)^m$.
The Hardy-Rellich inequality for $(-\Delta)^m$ was established by Owen \cite{O} who showed that if $\Omega$ is convex then
\begin{equation}
\int_{\Omega}(\Delta^{m/2}u)^2dx \geq A(m)\int_{\Omega}\frac{u^2}{d^{2m}}dx\; , \qquad u\in C^{\infty}_c(\Omega),
\label{owen}
\end{equation}
where
\[
A(m)=\frac{1^2\cdot 3^2\cdot\ldots\cdot  (2m-1)^2}{4^m}
\]
is sharp. Here and below we abuse the notation and write $\int(\Delta^{m/2}u)^2dx$ to stand for
$\int |\nabla\Delta^{(m-1)/2}u|^2dx$ when $m$ is odd. In the main theorems of this paper
we obtain two improvements of (\ref{owen}), a logarithmic and a geometric improvement.
To state our results, let us define the constants
\begin{eqnarray*}
&& B(m)=\frac{1}{4^m}\sum_{i=1}^m\prod_{\sscr{k=1}{k\neq i}}^m(2k-1)^2, \\
&&\Gamma(m)=\frac{N(N+2)\ldots (N+2m-2)}{1\cdot 3 \cdots (2m-1)}\Big(
\sum_{i=1}^m\frac{1}{2^{m+i}}\prod_{k=1}^i (2k-1)^2\Big)\alpha_N^{2m/N}.
\end{eqnarray*}
Our first theorem yields a logarithmic series improvement:
\begin{theorem}
Let $\Omega$ be convex and such that $d(x)$ is bounded in $\Omega$. Then there exists $D\geq sup_{\Omega}d(x)$ such that
\[
\int_{\Omega}(\Delta^{m/2}u)^2dx \geq A(m)\int_{\Omega}\frac{u^2}{d^{2m}}dx 
+B(m)\sum_{i=1}^{\infty}
 \int_{\Omega}\frac{u^2}{d^{2m}}X_1^2(d/D) \ldots X_i^2(d/D) dx \, , 
\]
for all functions $u\in C^{\infty}_c(\Omega)$.
\label{thm:1}
\end{theorem}
In the direction of geometric improvement we have
\begin{theorem}
Let $\Omega$ be bounded and convex. Then there holds
\[
\int_{\Omega}(\Delta^{m/2}u)^2dx \geq A(m)\int_{\Omega}\frac{u^2}{d^{2m}}dx
 +\Gamma(m)|\Omega|^{-2m/N}\int_{\Omega}u^2dx \, , 
\]
for all functions $u\in C^{\infty}_c(\Omega)$.
\label{thm:3}
\end{theorem}
For $m=2$ Theorem \ref{thm:1} recovers inequality (\ref{bt}), while for $m=1$ Theorem \ref{thm:3} recovers (\ref{hhl}).
The constant $B(m)$ of Theorem \ref{thm:1} is sharp; this is contained in the next theorem: we set
\[
 I_r[u]=\int_{\Omega}(\Delta^{m/2}u)^2dx - A(m)\int_{\Omega}\frac{u^2}{d^{2m}}dx 
-B(m)\sum_{i=1}^r \int_{\Omega}\frac{u^2}{d^{2m}}X_1^2(d/D) \ldots X_i^2(d/D) dx \, .
\]
\begin{theorem}
Let $r\geq 1$ and suppose that for some constants $C> 0$, $\theta\in\R$ and $D\geq \sup_{\Omega} d(x)$
the following inequality holds true,
\begin{equation}
I_{r-1}[u]   \geq  C \int_{\Omega}\frac{u^2}{d^{2m}}X_1^2(d/D)\ldots 
X_{r-1}^2(d/D)X_r^{\theta}(d/D)dx.
\label{qq}
\end{equation}
for all $u\in C^{\infty}_c(\Omega)$. Then $\ia$ $\theta\geq 2$. $\ib$ If $\theta = 2$ then $C \leq B(m)$.
\label{thm:2}
\end{theorem}
We point out that the value of $D$ does not affect the optimality of Theorem \ref{thm:3} since for any $D_1,D_2\geq\sup_{\Omega}d(x)$
there holds $\lim (X_i(d/D_1))/(X_i(d/D_2)) =1$ as $x\to\partial\Omega$.

Our proofs of Theorems \ref{thm:1} and \ref{thm:3} are surprisingly simple once some one-dimensional inequalities
are available. These inequalities are obtained in Section \ref{sec2}.
With these in hand the proof is completed using the mean-distance function introduced by Davies \cite{D},
as adapted in \cite{O}; this is carried out in Section \ref{sec3}.
What is significantly more involved is the proof of the optimality of the constant $B(m)$
in Theorem \ref{thm:2}. This is established in Section \ref{sec4}.

\section{One dimensional estimates} \label{sec2}

For $\gamma>-1$ we define the constants
\begin{eqnarray*}
&&A(m,\gamma)=\frac{(\gamma+1)^2(\gamma+3)^2\ldots (\gamma+2m-1)^2}{4^m}, \\
&&B(m,\gamma)=\frac{1}{4^m}\sum_{i=1}^m\prod_{\sscr{k=1}{k\neq i}}^m(\gamma+2k-1)^2\, , \\
&&\Gamma(m,\gamma)=\frac{N(N+2)\ldots (N+2m-2)}{(\gamma+1)(\gamma+3)\ldots (\gamma+2m-1)}
\biggl(\sum_{i=1}^m\frac{1}{2^{m+i}}\prod_{k=1}^i (\gamma+2k-1\Big)^2\biggl)\alpha^{2m/N}.
\end{eqnarray*}
Note that when $\gamma=0$ these reduce
to the constants $A(m)$, $B(m)$ and $\Gamma(m)$ defined in the introduction. In relation to the case $m=1$ of this definition,
throughout the paper we adopt the convention that empty sums equal zero and empty products equal one.

To simplify the notation we define
\begin{equation}
\zeta(s)= \sum_{i=1}^{\infty} X_1^2(s) \ldots X_i^2(s)\; , \quad s\in (0,1].
\label{oui}
\end{equation}
Throughout this section we fix an open interval $(0,2b)$
and let $\rho(t)=\min\{t,2b-t\}$, the distance of $t$
to the boundary of $\{0,2b\}$. We have
\begin{prop}
Let $m\geq 1$ be fixed. Then there exists $D\geq b$ such that for any $\gamma >-1$ and $\lambda\geq 0$ there holds
\begin{eqnarray}
&&\int_0^{2b}\!\!(1+\lambda\zeta(\rho/D)) \frac{(u^{(m)})^2}{\rho^{\gamma}(t)}dt \geq A(m,\gamma)
\int_0^{2b}\!\!\frac{u^2}{\rho^{\gamma +2m}}dt+\nonumber \\
&&\hspace{2cm}+\Big[B(m,\gamma)+\lambda A(m,\gamma)\Big]
\int_0^{2b}\!\!\frac{u^2}{\rho^{\gamma+2m}}\zeta(\rho/D)\, dt \, ,
\label{oned}
\end{eqnarray}
for all $u\in C^{\infty}_c(0,2b)$.
\label{p:oned}
\end{prop}
{\em Proof.} We use induction. For $m=1$ the result is contained in \cite[Theorem 1]{BT};
crucially, the constant $D$ does {\em not} depend on $\gamma$.
We assume that (\ref{oned}) is valid for $m-1$ (for the same $D$
and for {\em any} $\gamma>-1$) and writing for simplicity $\zeta$ for $\zeta(\rho(t)/D)$, we have
\begin{eqnarray*}
&&\qquad \int_0^{2b}(1+\lambda\zeta)\frac{(u^{(m)})^2}{\rho^{\gamma}}dt \\
&\geq& A(m-1,\gamma)\int_0^{2b}\frac{(u')^2}{\rho^{\gamma +2m-2}}dt +\\
&&+\Big[ B(m-1,\gamma)+\lambda A(m-1,\gamma)\Big]\int_0^{2b}\frac{(u')^2}{\rho^{\gamma +2m-2}}\zeta\, dt \\
&=&A(m-1,\gamma)\left\{ \int_0^{2b}\bigg( 1 +\Big[\lambda+\frac{B(m-1,\gamma)}{A(m-1,\gamma)}\Big]\zeta\bigg)
\frac{(u')^2}{\rho^{\gamma+2m-2}}dt\right\}\\
&\geq&A(m-1,\gamma)\left\{A(1,\gamma+2m-2)\int_0^{2b}\frac{u^2}{\rho^{\gamma+2m}}dt +\right.\\
&&\left.+\bigg[B(1,\gamma+2m-2)+\Big[\lambda+\frac{B(m-1,\gamma)}{A(m-1,\gamma)}\Big]
A(1,\gamma+2m-2)\bigg]
\int_0^{2b}\frac{u^2}{\rho^{\gamma+2m}}\zeta\, dt\right\} \\
&=&A(m-1,\gamma)A(1,\gamma+2m-2)\int_0^{2b}\frac{u^2}{\rho^{\gamma+2m}}dt +\\
&& +\left\{\ \Big[A(m-1,\gamma)B(1,\gamma+2m-2)+
B(m-1,\gamma)A(1,\gamma+2m-2)\Big] +\right. \\
&&+\lambda A(m-1,\gamma)A(1,\gamma+2m-2)
\int_0^{2b}\frac{u^2}{\rho^{\gamma+2m}}\zeta\, dt\, .
\end{eqnarray*}
Now, simple calculations together with the relations $A(1,\gamma)=(\gamma+1)^2/4$ and
$B(1,\gamma)=1/4$ show that
\begin{eqnarray*}
&&A(m,\gamma)=A(m-1,\gamma)A(1,\gamma+2m-2)\, , \\
&&B(m,\gamma)=A(m-1,\gamma)B(1,\gamma+2m-2)+B(m-1,\gamma)A(1,\gamma+2m-2)\, .
\end{eqnarray*}
This concludes the proof. $\hfill //$

\begin{lemma}
Let $\gamma>-1$ be fixed. Then
\begin{equation}
\int_0^{2b}\frac{(u')^2}{\rho^{\gamma}}dt\geq \frac{(\gamma+1)^2}{4}\int_0^{2b}\frac{u^2}{\rho^{\gamma+2}}dt
+\frac{(\gamma+1)^2}{4}\frac{1}{b^{\gamma+2}}\int_0^{2b}u^2dt  ,
\label{ei1}
\end{equation}
for all functions $u\in C^{\infty}_c(0,2b)$.
\label{lem:nw}
\end{lemma}
{\em Proof.} Let $u\in C^{\infty}_c(0,2b)$ be given and let $g$ be a continuous function on $(0,b)$. There holds
\begin{eqnarray*}
\int_0^{b}g'(\rho(t))u^2dt&=&g(b)u^2(b)-2\int_0^bg(\rho(t))uu'dt \\
&\leq&g(b)u^2(b) +\int_0^bg^2(\rho(t))\rho^{\gamma}u^2dt +\int_0^b\frac{(u')^2}{\rho^{\gamma}}dt,
\end{eqnarray*}
that is
\[ 
\int_0^b\frac{(u')^2}{\rho^{\gamma}}dt\geq \int_0^b\Big( g'(\rho(t))-g^2(\rho(t))\rho^{\gamma}\Big)u^2dt -g(b)u^2(b).
\]
Similarly,
\[ 
\int_b^{2b}\frac{(u')^2}{\rho^{\gamma}}dt\geq \int_b^{2b}\Big( g'(\rho(t))-g^2(\rho(t))\rho^{\gamma}\Big)u^2dt -g(b)u^2(b).
\]
Adding up we obtain
\[ 
\int_0^{2b}\frac{(u')^2}{\rho^{\gamma}}dt\geq \int_0^{2b}\Big( g'(\rho(t))-g^2(\rho(t))\rho^{\gamma}\Big)u^2dt -2g(b)u^2(b).
\]
Replacing $g(\cdot)$ by $g(\cdot)-g(b)$ we conclude that
\begin{equation}
\int_0^{2b}\frac{(u')^2}{\rho^{\gamma}}dt\geq \int_0^{2b}\Big( g'(\rho(t))-[g(\rho(t))-g(b)]^2\rho^{\gamma}\Big)u^2dt.
\label{har}
\end{equation}
Choosing
\[
g(s) =-\frac{\gamma+1}{2}s^{-\gamma-1},
\]
yields after some simple calculations
\begin{eqnarray}
\int_0^{2b}\frac{(u')^2}{\rho^{\gamma}}dt&\geq&\frac{(\gamma+1)^2}{4}\int_0^{2b}\frac{u^2}{\rho^{\gamma+2}}dt+
\frac{(\gamma+1)^2}{2}\int_0^{2b}\frac{u^2}{b^{\gamma+1}\rho}dt
-\frac{(\gamma+1)^2}{4}\int_0^{2b}\frac{\rho^{\gamma}u^2}{b^{2\gamma+2}}dt\nonumber \\
&\geq&\frac{(\gamma+1)^2}{4}\int_0^{2b}\frac{u^2}{\rho^{\gamma+2}}dt 
+\frac{(\gamma+1)^2}{4}\int_0^{2b}\frac{u^2}{b^{\gamma+2}}dt.
\label{ccc}
\end{eqnarray}
\vskip-3em\hfill //

\

For $\gamma>-1$ we define
\[
E(m,\gamma)=\sum_{i=1}^m\frac{1}{2^{m+i}}\prod_{k=1}^i (\gamma+2k-1)^2.
\]
\begin{prop}
For any $\gamma>-1$ there holds
\begin{equation}
\int_0^{2b}\frac{(u^{(m)})^2}{\rho^{\gamma}}dt\geq A(m,\gamma)\int_0^{2b}\frac{u^2}{\rho^{\gamma+2m}}dt
+E(m,\gamma)\frac{1}{b^{\gamma +2m}}\int_0^{2b}u^2dt,
\label{ele}
\end{equation}
for all functions $u\in C^{\infty}_c(0,2b)$.
\label{pr:br}
\end{prop}
{\em Proof.} For $m=1$ this has been proved in the last lemma. Assuming (\ref{ele}) to be true for $m-1$ we compute
\begin{eqnarray*}
\int_0^{2b}\frac{(u^{(m)})^2}{\rho^{\gamma}}dt&\geq&A(m-1,\gamma)\int_0^{2b}\frac{(u')^2}{\rho^{\gamma+2m-2}}dt
+E(m-1,\gamma)\int_0^{2b}\frac{(u')^2}{b^{\gamma+2m-2}}dt\\
&\geq&A(m-1,\gamma)\frac{(2m-1+\gamma)^2}{4}\int_0^{2b}\frac{u^2}{\rho^{\gamma+2m}}dt + \\
&&+
\bigg( A(m-1,\gamma)\frac{(2m-1+\gamma)^2}{4} +\frac{1}{2}\frac{E(m-1,\gamma)}{2}\bigg)\frac{1}{b^{\gamma+2m}}\int_0^{2b}u^2dt\, .
\end{eqnarray*}
The result follows if we note that
\[
A(m,\gamma)=A(m-1,\gamma)\frac{(2m-1+\gamma)^2}{4}\quad , \quad E(m,\gamma)=A(m,\gamma)+\frac{1}{2}E(m-1,\gamma).\]
\vskip-3em\hfill //

\

{\bf Remark.} We could use the intermediate inequality in (\ref{ccc}), hence obtaining $b^{-\gamma-1}\rho^{-1}$ instead of
$b^{-\gamma-2}$ in (\ref{ei1}). This would lead to a better constant $\hat{E}(m,\gamma)$, defined inductively by
\[
\hat{E}(1,\gamma)=\frac{(\gamma+1)^2}{4} \; , \qquad 
\hat{E}(m,\gamma)=\frac{(\gamma+1)^2}{4}\Big[ \hat{E}(m-1,\gamma+2) +\hat{E}(m-1,1)+A(m-1,1)
\Big].
\]

\section{Higher dimensions} \label{sec3}

Let $\Omega$ be a convex domain in $\R^N$.
We introduce some additional notation (see \cite{D,HHL}).
For $\omega\in S^{N-1}$ and $x\in\Omega$ we define
the following functions with values in $(0,+\infty]$:
\begin{eqnarray}
\tau_{\omega}(x)&=&\inf\{s>0 \; |\; x+s\omega\not\in\Omega\} \nonumber \\
\rho_{\omega}(x)&=&\min\{\tau_{\omega}(x),\tau_{-\omega}(x)\}\\
b_{\omega}(x)&=&\frac{1}{2}(\tau_{\omega}(x)+\tau_{-\omega}(x)).\label{ari}
\nonumber
\end{eqnarray}

We can now prove Theorems \ref{thm:1} and \ref{thm:3}.

{\em Proof of Theorem \ref{thm:1}.} Let $u\in C^{\infty}_c(\Omega)$ be given.
Let us fix a direction $\omega\in S^{N-1}$ and let $\Omega_{\omega}$
be the orthogonal projection of $\Omega$ on the hyperplane perpendicular to
$\omega$. For each $z\in\Omega_{\omega}$ we apply Proposition \ref{p:oned} (with $\gamma=0$) on
the segment defined by $z$ and $\omega$. By continuity and compactness, $D$ can be chosen to be
independent of $\omega$.
We then integrate over $z\in\Omega_{\omega}$ and using the convexity of $\Omega$ we conclude that
\[
\int_{\Omega}(\partial_{\omega}^mu)^2dx\geq
A(m)\int_{\Omega}\frac{u^2}{\rho_{\omega}^{2m}}dx
+B(m)\int_{\Omega}\frac{u^2}{\rho_{\omega}^{2m}}\zeta(\rho_{\omega}(x)/D)dx\, .
\]
Since $\zeta$ is an increasing function, this implies
\begin{equation}
\int_{\Omega}(\partial_{\omega}^mu)^2dx\geq
A(m)\int_{\Omega}\frac{u^2}{\rho_{\omega}^{2m}}dx 
+B(m)\int_{\Omega}\frac{u^2}{\rho_{\omega}^{2m}}\zeta(d(x)/D)dx\, .
\label{thi}
\end{equation}
We now integrate over $\omega\in S^{N-1}$. It is shown in
\cite{O} that
\begin{equation}
\int_{S^{N-1}}\int_{\Omega}(\partial_{\omega}^mu)^2dx\, dS(\omega)
=C(m,N)\int_{\Omega}(\Delta^{m/2}u)^2dx\; ,
\label{o1}
\end{equation}
where
\[
C(m,N)=\frac{1\cdot 3\ldots (2m-1)}{N(N+2)\ldots (N+2m-2)}.
\]
In the same article it was shown that the convexity of $\Omega$ implies
\begin{equation}
\int_{S^{N-1}}\frac{dS(\omega)}{\rho_{\omega}^{2m}(x)}\geq C(m,N)\frac{1}{d^{2m}(x)}.
\label{o2}
\end{equation}
Combining (\ref{thi}), (\ref{o1}) and (\ref{o2}) we obtain the stated inequality.
$\hfill //$

\

{\bf Proof of Theorem \ref{thm:3}.} Let $u\in C^{\infty}_c(\Omega)$ be given.
Arguing as before, but using now Proposition \ref{pr:br} instead of Proposition \ref{p:oned}, we have
\[
\int_{\Omega}(\partial_{\omega}^mu)^2dx\geq
A(m)\int_{\Omega}\frac{u^2}{\rho_{\omega}^{2m}}dx 
+E(m)\int_{\Omega}\frac{u^2}{b_{\omega}^{2m}}dx\, , \qquad \omega\in S^{N-1}.
\]
Integrating over $\omega\in S^{N-1}$ and using (\ref{o1}) and (\ref{o2}) yields
\begin{equation}
\int_{\Omega}(\Delta^{m/2}u)^2dx \geq A(m)\int_{\Omega}\frac{u^2}{d^{2m}}dx 
+\frac{E(m)}{C(m,N)}\int_{\Omega}\int_{S^{N-1}}\frac{u^2}{b_{\omega}^{2m}}dS(\omega)dx\, .
\label{thiv}
\end{equation}
But \cite[Lemma 2.1]{T} the convexity of $\Omega$ implies that
\begin{equation}
\int_{S^{N-1}}\frac{1}{b_{\omega}(x)^{2m}}dS(\omega)\geq \Big( \frac{|\Omega|}{a_N}\Big)^{-2m/N}.
\label{thiv1}
\end{equation}
Combining (\ref{thiv}) and (\ref{thiv1}) and observing that
\[
\Gamma(m)=\frac{E(m)}{C(m,N)}\alpha_N^{2m/N},
\]
concludes the proof of the theorem. $\hfill //$

%%%%%%%%%%%%%%%%%%%%%%%%%%%%%%%%%%%%%%

\section{Optimality of the constants} \label{sec4}

This section is considerably more technical than the previous ones. Our main purpose will be the
computation of $I_{r-1}[u]$ for an appropriate test function $u$.
Throughout the section we shall repeatedly use the differentiation rule
\begin{equation}
 \frac{d}{dt} X_i^{\beta}(t)=\frac{\beta}{t}X_1(t)X_2(t)\ldots X_{i-1}(t)X_i^{1+\beta}(t),\qquad i=1,2,\ldots,\quad \beta\in\R,
\label{diffrule}
\end{equation}
which is easily proved by induction.

Let $m\in\N$. We recall our convention about empty sums or products and define the functions
\begin{eqnarray*}
&& \sigma^{(m)}_0(x)=x(x-1)\ldots (x-m+1)\;\; , \;\; \sigma^{(m)}_1(x)=\sum_{i=1}^m\prod_{k\neq i}(x-k+1) \\
&& \sigma^{(m)}_2(x)=\sum_{1\leq i<j\leq r}^m\prod_{k\neq i,j}(x-k+1) .
\end{eqnarray*}

\begin{lemma}
Let $s_0,s_1,\ldots,s_r\in\R$ and $u(t)=t^{s_0}X_1^{s_1}\ldots X_r^{s_r}$. Let
\[
Y_{ij}=X_1^2\ldots X_i^2X_{i+1}\ldots X_j\, \qquad 0\leq i\leq j\leq r,
\]
with the conventions $Y_{00}=1$, $Y_{ii}=X_1^2\ldots X_i^2$, $Y_{0j}=X_1\ldots X_j$.
Then there holds
\begin{equation}
u^{(m)}(t)=t^{s_0-m}X_1^{s_1}\ldots X_r^{s_r}
\sum_{0\leq i\leq j\leq r}c_{ij}^{(m)}Y_{ij}(t) + t^{s_0-m}O(X_1^{s_1+3}X_2^{s_2}\ldots X_r^{s_r}),
\label{mder}
\end{equation}
where:
\[
\begin{array}{lll}
 c_{00}^{(m)}=\sigma^{(m)}_0(s_0),  & c_{0j}^{(m)}=s_j\sigma^{(m)}_1(s_0)  ,  & j\geq 1, \\
 c_{ii}^{(m)}=s_i(s_i+1)\sigma^{(m)}_2(s_0) \; , \quad  1\leq i\leq r, & 
c_{ij}^{(m)}=(2s_i+1)s_j\sigma^{(m)}_2(s_0) , &  1\leq i<j\leq r\, .
\end{array}
\]
\label{lem:mder}
\end{lemma}
{\em Proof.} We use induction. When $m=1$ (\ref{mder}) follows directly from (\ref{diffrule}).
We assume that
\[
u^{(m-1)}(t)=t^{s_0-m+1}X_1^{s_1}\ldots X_r^{s_r}
\sum_{0\leq i\leq j\leq r}c_{ij}^{(m-1)}Y_{ij}(t) + t^{s_0-m+1}O(X_1^{s_1+3}X_2^{s_2}\ldots X_r^{s_r}),
\]
We differentiating and again use (\ref{diffrule}). The $t^{s_0-m+1}O(X_1^{s_1+3}X_2^{s_2}\ldots X_r^{s_r})$
will give a term $t^{s_0-m}O(X_1^{s_1+3}X_2^{s_2}\ldots X_r^{s_r})$.
After some simple calculations we obtain modulo $O(X_1^{s_1+3}X_2^{s_2}\ldots X_r^{s_r})$,
\begin{eqnarray*}
u^{(m)}(t)&=&t^{s_0-m}X_1^{s_1}\ldots X_r^{s_r}\Bigg\{
\sum_{0\leq i\leq j\leq r}c_{ij}^{(m-1)}(s_0-m+1)Y_{ij} +\sum_{j=i}^rc_{00}^{(m-1)}s_jY_{0j}\\
&&+\sum_{j=1}^r\sum_{k=1}^rc_{0j}^{(m-1)}Y_{kj} +\sum_{j=1}^r\sum_{k=j+1}^rc_{0j}^{(m-1)}s_kY_{jk}\Bigg\} \\
&=&t^{s_0-m}X_1^{s_1}\ldots X_r^{s_r}\Bigg\{
(s_0-m+1)c_{00}^{(m-1)}Y_{00} \\
&&+ \sum_{j=1}^r \Big[ (s_0-m+1)c_{0j}^{(m-1)} +s_jc_{00}^{(m-1)}\Big]Y_{0j} \\
&& +\sum_{i=1}^r \Big[(s_0-m+1)c_{ii}^{(m-1)} +(s_i+1)c_{0i}^{(m-1)} \Big] Y_{ii} \\
&&+\sum_{1\leq i<j\leq r}\Big[ (s_0-m+1)c_{ij}^{(m-1)} +(s_i+1)c_{0j}^{(m-1)}+s_jc_{0i}^{(m-1)}\Big]Y_{ij}\Bigg\}.
\end{eqnarray*}
The proof is concluded by observing that the constants $c_{ij}^{(k)}$, $0\leq i\leq j\leq r$,
satisfy the induction relations
\[
\begin{array}{lll}
& c_{00}^{(m)}=(s_0-m+1)c_{00}^{(m-1)},  & \\[0.2cm]
& c_{0j}^{(m)}=(s_0-m+1)c_{0j}^{(m-1)} +s_jc_{00}^{(m-1)}, &  1\leq j\leq r,\\[0.2cm]
& c_{ii}^{(m)}=(s_0-m+1)c_{ii}^{(m-1)} +(s_i+1)c_{0i}^{(m-1)}, & 1\leq i\leq r,\\[0.2cm]
& c_{ij}^{(m)}=(s_0-m+1)c_{ij}^{(m-1)} +(s_i+1)c_{0j}^{(m-1)} +s_jc_{0i}^{(m-1)}, &  1\leq i< j\leq r. 
\end{array}
\]
$\hfill //$

In the sequel we shall denote the constants $c_{ij}^{(m)}$ simply by $c_{ij}$, since only
the $m$th order derivative of $u$ will appear.
Similarly, we shall write $\sigma_i(x)$ instead of $\sigma_i^{(m)}(x)$, $i=0,1,2$.
Let $s_0>(2m-1)/2$, $s_1,\ldots,s_r\in\R$ be fixed.
For $0\leq i\leq j\leq r$ we define
\begin{eqnarray*}
\Gamma_{ij}&=&\int_0^1 t^{2s_0-2m}X_1^{2s_1}\ldots X_r^{2s_r}Y_{ij}dt \\
&=&\int_0^1 t^{2s_0-2m}X_1^{2s_1+2}\!\!\!\!\ldots
X_i^{2s_i+2}X_{i+1}^{2s_{i+1}+1}\!\!\ldots X_j^{2s_j+1}
X_{j+1}^{2s_{j+1}}\!\!\!\!\ldots X_r^{2s_r} dt.
\end{eqnarray*}
\begin{lemma}
Let $u(t)=t^{s_0}X_1^{s_1}\ldots X_r^{s_r}$. There holds
\begin{equation}
I_{r-1}[u]=\sum_{0\leq i\leq j\leq r}a_{ij}\Gamma_{ij} + \int_0^1 t^{2s_0-2m}O(X_1^{s_1+3}X_2^{s_2}\ldots X_r^{s_r})dt.
\label{rain}
\end{equation}
where
\begin{equation}
\begin{array}{llll}
a_{00} =c_{00}^2-\alpha(m) ,&   &  a_{0j} =2c_{00}c_{0j}, & 1\leq j\leq r,  \\[0.2cm]
a_{ii}= c_{0i}^2+2c_{00}c_{ii} -\beta(m) , &  1\leq i\leq r-1, &a_{rr}=c_{0r}^2+2c_{00}c_{rr} ,\\[0.2cm]
a_{ij} =2c_{00}c_{ij}+2c_{0i}c_{0j}, & 1\leq i<j\leq r.
\end{array}
\label{aij}
\end{equation}
\label{trian}
\end{lemma}
{\em Proof.} From Lemma \ref{lem:mder} we have modulo
$\int_0^1 t^{2s_0-2m}O(X_1^{s_1+3}X_2^{s_2}\ldots X_r^{s_r})dt$,
\[
 \int_0^1 (u^{(m)})^2dt = \int_0^1 t^{2s_0-2m}X_1^{2s_1}\ldots X_r^{2s_r}
\bigg( \sum_{0\leq i\leq j\leq r}c_{ij}Y_{ij}\bigg)^2 dt \, .
\]
We expand the square and hence obtain a linear combination of terms of the form
$\int_0^1 t^{2s_0-2m}X_1^{2s_1}\ldots X_r^{2s_r}Y_{ij}Y_{kl}dt$, where
$0\leq i\leq j\leq r$, $0\leq k\leq l\leq r$. Now, we observe that
$Y_{ij}Y_{kl}=O(X_1^3)$
unless (1) $i=j=0$ or (2) $k=l=0$ or (3) $i=k=0$. Hence, denoting by $S$ the last parenthesis above we have
\begin{eqnarray*}
S&=&c_{00}^2 +2\!\!\!\!\!\sum_{\sscr{0\leq i\leq j\leq r}{(i,j)\neq (0,0)}}\!\!\!\!\! c_{00}c_{ij}Y_{ij} 
+\sum_{j,l=1}^rc_{0j}c_{0l}Y_{0j}Y_{0l} + O(X_1^3) \\
&=&c_{00}^2 +2\sum_{j=1}^rc_{00}c_{0j}Y_{0j}  +2\!\!\!\!\!\sum_{1\leq i\leq j\leq r}\!\!\!\!\! c_{00}c_{ij}Y_{ij}
+\sum_{i=1}^rc_{0i}^2Y_{0i}^2+2\!\!\!\!\!\sum_{1\leq i<j\leq r}\!\!\!\!\! c_{0i}c_{0j}Y_{0i}Y_{0j} + O(X_1^3).
\end{eqnarray*}
Using the fact that $Y_{0i}Y_{0j}=Y_{ij}$, $i\leq j$, we thus conclude that
\[
S=c_{00}^2 +2\sum_{j=1}^rc_{00}c_{0j}Y_{0j} +\sum_{i=1}^r(c_{0i}^2+2c_{00}c_{ii})Y_{ii} 
+2\!\!\!\!\!\sum_{1\leq i<j\leq r}\!\!\!\!\! (c_{00}c_{ij}+2c_{0i}c_{0j})Y_{ij}.
\]
The proof is complete if we recall that
\[  \int_0^1 \frac{u^2}{t^{2m}}dt= \Gamma_{00}
\mbox{ and } \int_0^1 \frac{u^2}{t^{2m}}X_1^2\ldots X_i^2 dt =\Gamma_{ii}\; , \quad 1\leq i\leq r-1\, .
\]
\vskip-3em\hfill //

\

Up to this point the parameters $s_0,s_1,\ldots, s_r$ where arbitrary subject only to
$s_0>(2m-1)/2$. We now make a more specific choice, taking
\[
s_0=\frac{2m-1+\epsilon_0}{2}\; , \qquad s_j=\frac{-1+\epsilon_j}{2}\; , \quad 1\leq j\leq r\, ,
\]
where $\epsilon_0,\ldots,\epsilon_r$ are small parameters.
We consider the functional $I_{r-1}[u]$ as a function of these parameters and intend
to take succesively the limits $\epsilon_0\searrow 0,\ldots ,\epsilon_r\searrow 0$. In taking these
limits we shall ignore terms that are bounded uniformly in the $\epsilon_i$'s. In order to distinguish
such terms we shall make use of the following fact: we have \cite[(3.8)]{BFT2}:
\begin{equation}
\int_0^1  t^{-1+\epsilon_0}X_1^{1+\epsilon_1}\ldots X_r^{1+\epsilon_r}dt <\infty 
\Longleftrightarrow\; \left\{
\begin{array}{ll}
& \epsilon_0>0 \\
\mbox{ or}&\mbox{$\epsilon_0=0$ and $\epsilon_1>0$}\\
\mbox{ or}&\mbox{$\epsilon_0=\epsilon_1=0$ and $\epsilon_2>0$}\\
& \cdots \\
\mbox{ or}&\mbox{$\epsilon_0=\epsilon_1=\ldots =\epsilon_{r-1}=0$ and $\epsilon_r>0$.}
\end{array}\right.
\label{finite}
\end{equation}

For the terms that diverge as the $\epsilon_i$'s tend to zero,
we shall need some quantitive information on the rate of divergence. This is contained in
the following
\begin{lemma}
For any $\beta<1$ there exists $c_{\beta}>0$ such that
\begin{eqnarray*}
&\ia& \int_0^1 t^{-1+\epsilon_0}X_1^{\beta}dt
\leq c_{\beta}\epsilon_0^{-1+\beta}\; ,\\
&\ib & \int_0^1 t^{-1}X_1\ldots X_{i-1}X_i^{1+\epsilon_i}X_{i+1}^{\beta}dt 
\leq c_{\beta}\epsilon_i^{-1+\beta}\; , \qquad 1\leq i\leq r-1\, .
\end{eqnarray*}
\label{diverge}
\end{lemma}
{\em Proof.} $\ia$ Setting $s=\epsilon_0^{-1}X_1(t)$ we have $t=\exp(1-\epsilon_0^{-1}s^{-1})$,
$ds=\epsilon_0^{-1}t^{-1}X_1^2dt$, and therefore
\begin{eqnarray*}
\int_0^1 t^{-1+\epsilon_0}X_1^{\beta}dt
&=&e^{\epsilon_0}\epsilon_0^{-1+\beta}\int_0^{\frac{1}{\epsilon_0}}e^{-\frac{1}{s}}s^{-2+\beta}ds\\
&\leq& e^{\epsilon_0}\epsilon_0^{-1+\beta}\int_0^{\infty}e^{-\frac{1}{s}}s^{-2+\beta}ds.
\end{eqnarray*}
$\ib$ Similarly, we set $s=\epsilon_i^{-1}X_{i+1}(t)$. Then
\[
X_i(t)=\exp(1-\epsilon_i^{-1}s^{-1})\;\; , \qquad  ds=\epsilon_i^{-1}t^{-1}X_1\ldots X_iX_{i+1}^2dt.
\]
Hence (\ref{diffrule}) gives
\[
\int_0^1 t^{-1}X_1\ldots X_{i-1}X_i^{1+\epsilon_i}X_{i+1}^{\beta}dt=
e^{\epsilon_i}\epsilon_i^{-1+\beta}\int_0^{\frac{1}{\epsilon_i}}e^{-\frac{1}{s}}s^{-2+\beta}ds,
\]
yielding the stated estimate. $\hfill //$

We shall also need the following
\begin{lemma}
$\ia$ There holds
\[
 \epsilon_0^2\Gamma_{00} -2\epsilon_0\sum_{j=i+1}^r(1-\epsilon_j)\Gamma_{0j} =
\sum_{i=1}^r(\epsilon_i-\epsilon_i^2)\Gamma_{ii} -\sum_{1\leq i<j\leq r}(1-\epsilon_j)(1-2\epsilon_i)\Gamma_{ij}
+O(1),
\]
where the $O(1)$ is uniform in $\epsilon_0,\ldots,\epsilon_r$. \nl
$\ib$
Let $i\geq 0$ and (if $i\geq 1$) assume that $\epsilon_0=\ldots=\epsilon_{i-1}=0$. Then
\[
\epsilon_i\Gamma_{ii} = \sum_{j=i+1}^r(1-\epsilon_j)\Gamma_{ij} +O(1),
\]
where the O(1) is uniform in $\epsilon_i,\ldots,\epsilon_r$.
\label{lem:vil}
\end{lemma}
{\em Proof.} The two parts of the lemma have been proved in
\cite[p184]{BFT2} and \cite[p181]{BFT2} respectively. $\hfill //$

{\bf Remark.} We are now in position to prove Theorem \ref{thm:2}, but before proceeding
some comments are necessary.
The proof of the theorem is local: we fix a point $x_0\in\partial\Omega$ and work
entirely in a small ball $B(x_0,\delta)$ using a cut-off function $\phi$. The sequence of functions
that is used is then given by
\[
u(x)=\phi(x) d(x)^{\frac{-1+2m+\epsilon_0}{2}}X_1(d(x)/D)^{\frac{-1+\epsilon_1}{2}}\ldots
X_r(d(x)/D)^{\frac{-1+\epsilon_r}{2}}, \qquad (\epsilon_0,\ldots,\epsilon_r>0)
\]
and, as already mentioned, we take the successive limits $\epsilon_0\searrow 0,\ldots,\epsilon_r\searrow 0$;
in taking this limits, we work modulo terms that are bounded uniformly in the remaining $\epsilon_i$'s.
Such are any terms that contain derivatives of $\phi$; such are also any terms that contain second-order
derivatives of $d(x)$. Such terms involve necessarily $\Delta d$ and are dealt with using the fact that
$d\Delta d=O(d)$ as $x\to\partial\Omega$; this prevents the appearence
of any derivatives of $d(x)$ of order higher than two and so no such information is needed. These
considerations are to a large extent the justification of the fact that,
for the proof of Theorem \ref{thm:2} we can, without any loss of generality, restrict ourselves to the one-dimensional case.
We shall thus take $\Omega=(0,1)$, and consider the sequence
\[
u(t)=t^{\frac{-1+2m+\epsilon_0}{2}}X_1(t)^{\frac{-1+\epsilon_1}{2}}\ldots
X_r(t)^{\frac{-1+\epsilon_r}{2}}.
\]
discussed earlier; multiplication by an appropriate cut-off function shows that $u$ lies in the
appropriate Sobolev space.
Note that $u$ does not vanish at $t=1$, but the cut-off function $\phi$ would take care of that.
For a complete picture of what the full proof would look like, we refer to \cite{BT} where
the case $m=2$ has been carried out in every detail.

\

{\em Proof of Theorem \ref{thm:2}} (see also the remark above) We define
\begin{equation}
u(t)=t^{\frac{-1+2m+\epsilon_0}{2}}X_1(t)^{\frac{-1+\epsilon_1}{2}}\ldots
X_r(t)^{\frac{-1+\epsilon_r}{2}},
\label{u}
\end{equation}
where $\epsilon_0,\ldots,\epsilon_r$ are small positive parameters.
For the reader's convenience we recall from Lemma \ref{trian} that
\begin{equation}
I_{r-1}[u]=\sum_{0\leq i\leq j\leq r}a_{ij}\Gamma_{ij} +O(1),
\label{ma}
\end{equation}
where the $O(1)$ is uniform in $\epsilon_0,\ldots\epsilon_r$ (by (\ref{finite})) and
the constants $a_{ij}$ are given by
\begin{equation}
\begin{array}{llll}
a_{00} =c_{00}^2-A(m) ,&   &  a_{0j} =2c_{00}c_{0j}, & 1\leq j\leq r,  \\[0.2cm]
a_{ii}= c_{0i}^2+2c_{00}c_{ii} -B(m) , &  1\leq i\leq r-1, &a_{rr}=c_{0r}^2+2c_{00}c_{rr} ,\\[0.2cm]
a_{ij} =2c_{00}c_{ij}+2c_{0i}c_{0j}, & 1\leq i<j\leq r.
\end{array}
\label{aij1}
\end{equation}
The $c_{ij}$'s are given by
\[
\begin{array}{lll}
 c_{00}=\sigma_0(s_0),  & c_{0j}=s_j\sigma_1(s_0)  ,  & j\geq 1, \\
 c_{ii}=s_i(s_i+1)\sigma_2(s_0) \; , \quad  1\leq i\leq r, & 
c_{ij}=(2s_i+1)s_j\sigma_2(s_0) , &  1\leq i<j\leq r\, .
\end{array}
\]
where, in turn,
\[
s_0=\frac{2m-1+\epsilon_0}{2}\; , \qquad s_j=\frac{\epsilon_j-1}{2}\; , \quad 1\leq j\leq r,
\]
and
\begin{eqnarray*}
&& \sigma_0(x)=x(x-1)\ldots (x-m+1)\;\; , \;\; \sigma_1(x)=\sum_{i=1}^m\prod_{k\neq i}(x-k+1) \\
&& \sigma_2(x)=\sum_{1\leq i<j\leq m}\prod_{k\neq i,j}(x-k+1) .
\end{eqnarray*}
We observe that
\[
\sigma_0'(x)=\sigma_1(x)\quad ,  \qquad \sigma_1'(x)=2\sigma_2(x).
\]
We now let $\epsilon_0\searrow 0$ in (\ref{ma}). It follows from (\ref{finite}) that all $\Gamma_{ij}$'s
with $i\geq 1$ have finite limits.
As for the remaining terms $\Gamma_{0j}$\, , applying Lemma \ref{diverge} with $\beta=-3/2$ (for $j=0$)
and with $\beta=-1/2$ (for $j\geq 1$) we obtain respectively
\begin{eqnarray}
\Gamma_{00}&=&\int_0^1t^{-1+\epsilon_0}X_1^{-1+\epsilon_1}\ldots\ldots X_r^{-1+\epsilon_r}dt \nonumber\\
&\leq &c\int_0^1 t^{-1+\epsilon_0}X_1^{-\frac{3}{2}}dt \label{j00}\\
&\leq &c\epsilon_0^{-\frac{5}{2}} \nonumber
\end{eqnarray}
and
\begin{eqnarray}
\Gamma_{0j}&=&\int_0^1t^{-1+\epsilon_0}X_1^{\epsilon_1}\ldots X_j^{\epsilon_j}X_{j+1}^{-1+\epsilon_{j+1}}
\ldots X_r^{-1+\epsilon_r}dt \nonumber\\
&\leq &c\int_0^1 t^{-1+\epsilon_0}X_1^{-\frac{1}{2}}dt \label{j0}\\
&\leq &c\epsilon_0^{-\frac{3}{2}},\nonumber
\end{eqnarray}
where in both cases $c>0$ is independent of $\epsilon_1,\ldots,\epsilon_r$.
Now, we think of the contants $a_{0j}$ and $c_{0j}$ as functions of $\epsilon_0$, writting
$a_{0j}=a_{0j}(\epsilon_0)$, $c_{0j}=c_{0j}(\epsilon_0)$
and considering $\epsilon_1,\ldots,\epsilon_r$ as small positive parameters.
Using Taylor's theorem we shall expand the coefficient $a_{0j}$ of $\Gamma_{0j}$, $j= 0$ (resp. $j\geq 1$) in powers
of $\epsilon_0$, and relation (\ref{j00}) (resp. (\ref{j0})) shows that we can discard powers
with exponent $\geq 3$ (resp. $\geq 2$). We compute the remaining ones. Denoting by $A_{k,0j}$ the
coefficient of $\epsilon_0^k$ in $a_{0j}$ we have:

{\em - Constant term in $a_{00}$:} We have $A_{0,00}=a_{00}(0)=c_{00}^2(0)-A(m)=0$.

{\em - Coefficient of $\epsilon_0$ in $a_{00}$:} We have $c_{00}(\epsilon_0)=\sigma_0(\frac{2m-1+\epsilon_0}{2})$
and therefore $c_{00}'(0)=\frac{1}{2}\sigma_1(\frac{2m-1}{2})$. Hence
$a_{00}'(\epsilon_0)=2c_{00}(\epsilon_0)c_{00}'(\epsilon_0)=
\sigma_0(\frac{2m-1+\epsilon_0}{2})\sigma_1(\frac{2m-1+\epsilon_0}{2})$ and the coefficient is
\[ A_{1,00}=a_{00}'(0)=\sigma_0(\frac{2m-1}{2})\sigma_1(\frac{2m-1}{2}).\]
We henceforth write $\sigma_i$ for $\sigma_i((2m-1)/2)$, $i=0,1,2$.

{\em - Coefficient of $\epsilon_0^2$ in $a_{00}$:} The coefficient is
\[ A_{2,00}=\frac{1}{2}a_{00}''(0)=[c_{00}'(0)]^2+c_{00}(0)c_{00}''(0)=\frac{1}{4}\sigma_1^2
+\frac{1}{2}\sigma_0\sigma_2.
\]
{\em - Constant term in $a_{0j}$, $j\geq 1$:} This is
\[
 A_{0,0j}=a_{0j}(0)=2c_{00}(0)c_{0j}(0)=-(1-\epsilon_j)\sigma_0 \sigma_1\, .
\]
{\em - Coefficient of $\epsilon_0$ in $a_{0j}$:} This is
\begin{eqnarray*}
A_{1,0j}&=&a_{0j}'(0)=2c_{00}'(0)c_{0j}(0)+2c_{00}(0)c_{0j}'(0)\\
&=&-\frac{1}{2}(1-\epsilon_j)\sigma_1^2
-(1-\epsilon_j)\sigma_0\sigma_2.
\end{eqnarray*}
Now, we observe that $A_{0,0j}=-(1-\epsilon_j)A_{1,00}$. Hence $\ib$ of Lemma \ref{lem:vil} implies that
\begin{equation}
A_{1,00}\epsilon_0\Gamma_{00} +\sum_{j=1}^r A_{0,0j}\Gamma_{0j} =O(1)
\label{p1}
\end{equation}
uniformly in $\epsilon_1,\ldots,\epsilon_r$. Similarly, we observe that $A_{1,0j}=-2(1-\epsilon_j)A_{2,00}$.
Hence, by $\ia$ of Lemma \ref{lem:vil}, the remaining `bad' terms when combined give
\begin{eqnarray}
&&\hspace{-2cm}A_{2,00}\epsilon_0^2\Gamma_{00} +\epsilon_0\sum_{j=1}^r A_{1,0j}\Gamma_{0j}=\nonumber \\
&=&A_{2,00}\Big( \epsilon_0^2\Gamma_{00} -2\epsilon_0\sum_{j=1}^r (1-\epsilon_j)\Gamma_{0j}\Big)\label{hh}\\
&=&A_{2,00}\Big(
\sum_{i=1}^r(\epsilon_i-\epsilon_i^2)\Gamma_{ii} -\sum_{1\leq i<j\leq r}(1-\epsilon_j)(1-2\epsilon_i)\Gamma_{ij}
\Big)+O(1),
\nonumber
\end{eqnarray}
uniformly in $\epsilon_1,\ldots,\epsilon_r$. Note that the right-hand side of (\ref{hh}) has a finite limit as
$\epsilon_0\searrow 0$.
Combining (\ref{aij1}) , (\ref{p1}) and (\ref{hh}) we conclude that,
after letting $\epsilon_0\searrow 0$, we are left with
\begin{eqnarray}
&&\qquad I_{r-1}[u] \nonumber \\
&=&\sum_{i=1}^r\bigg( a_{ii}+ A_{2,00}(\epsilon_i-\epsilon_i^2)\bigg)\Gamma_{ii}
+\sum_{1\leq i< j\leq r}\bigg( a_{ij} -A_{2,00}(1-\epsilon_j)(1-2\epsilon_i)  \bigg)\Gamma_{ij}
 +O(1) \nonumber \\
&=:&\sum_{i=1}^rb_{ii}\Gamma_{ii} + \sum_{1\leq i<j\leq r}b_{ij}\Gamma_{ij}+O(1)\; , \quad\quad (\epsilon_0=0),
\label{alv}
\end{eqnarray}
where the $O(1)$ is uniform in $\epsilon_1,\ldots,\epsilon_r$.

We next let $\epsilon_1\searrow 0$ in (\ref{alv}). It follows from (\ref{finite}) that all the $\Gamma_{ij}$'s have
a finite limit, except those with $i=1$ which diverge to $+\infty$. The latter terms are again estimated with the aid
of Lemma \ref{diverge}, this time with $i=1$. Part $\ia$ of the lemma (with $\beta=-3/2$) yields
\begin{eqnarray}
\Gamma_{11}&=&\int_0^1t^{-1}X_1^{1+\epsilon_1}X_2^{-1+\epsilon_2}\ldots X_r^{-1+\epsilon_r}dt \nonumber\\
&\leq &c\int_0^1 t^{-1}X_1^{1+\epsilon_1}X_2^{-\frac{3}{2}}dt \label{j10}\\
&\leq &c\epsilon_1^{-\frac{5}{2}}\, , \nonumber
\end{eqnarray}
uniformly in $\epsilon_2,\ldots,\epsilon_r$.
For $j\geq 2$ it also yields (now with $\beta=-1/2$)
\begin{eqnarray}
\Gamma_{1j}&=&\int_0^1t^{-1}X_1^{1+\epsilon_1}X_2^{\epsilon_2}
\ldots X_j^{\epsilon_j}X_{j+1}^{-1+\epsilon_{j+1}}\ldots X_r^{-1+\epsilon_r}dt \nonumber\\
&\leq &c\int_0^1 t^{-1}X_1^{1+\epsilon_1}X_2^{-\frac{1}{2}}dt \label{j1}\\
&\leq &c\epsilon_1^{-\frac{3}{2}},\nonumber
\end{eqnarray}
again, uniformly in $\epsilon_2,\ldots,\epsilon_r$.
We think of the coefficients $b_{1j}$ and $a_{1j}$ as functions of $\epsilon_1$
and we expand these in powers of $\epsilon_1$.
Estimate (\ref{j10}) (resp. (\ref{j1})) implies that only the terms $1,\epsilon_1$ and $\epsilon_1^2$
(resp. 1 and $\epsilon_1$) give contributions for $\Gamma_{11}$
(resp. $\Gamma_{1j}$, $j\geq 2$) that do not vanish as $\epsilon_1\searrow 0$.
We shall compute the coefficients of these terms; note that $c_{00}$ is now treated simply as a constant.
Denoting by $B_{k,1j}$ the coefficient of $\epsilon_1^k$ in $b_{1j}$, $j\geq 1$, we have:

{\em - Constant term in $b_{11}$:} For $\epsilon_1=0$ we have $s_1=-1/2$. Hence
\begin{eqnarray*}
B_{0,11}&=&b_{11}(0)\\
&=&c_{01}^2(0)+2c_{00}c_{11}(0) -B(m)\\
&=&\frac{1}{4}\sigma_1^2 -\frac{1}{2}\sigma_0\sigma_2 - B(m) \\
&=&\frac{1}{4}\Big(\sum_{i=1}^m\prod_{k\neq i}\frac{2k-1}{2}\Big)^2
-\frac{1}{2}\Big(\prod_{k=1}^m\frac{2k-1}{2}\Big)
\Big(\sum_{1\leq i<j\leq m}\prod_{k\neq i,j}\frac{2k-1}{2}\Big) -\\
&&-\frac{1}{4}\sum_{i=1}^m \prod_{k\neq i}\Big(\frac{2k-1}{2}\Big)^2 .
\end{eqnarray*}
This is zero as is seen by expanding the square:
\[\Big(\sum_{i=1}^m\prod_{k\neq i}\frac{2k-1}{2}\Big)^2 
=\sum_{i=1}^m\prod_{k\neq i}\Big(\frac{2k-1}{2}\Big)^2 
+2\sum_{i<j}\Big(\prod_{k\neq i}\frac{2k-1}{2}\Big)\Big(\prod_{k\neq j}\frac{2k-1}{2}\Big).
\]
{\em - Coefficient of $\epsilon_1$ in $b_{11}$:} We have
\begin{eqnarray*}
b_{11}'(\epsilon_1)&=&a_{11}'(\epsilon_1)+A_{2,00}-2A_{2,00}\epsilon_1\\
&=&2c_{01}(\epsilon_1)c_{01}'(\epsilon_1)+2c_{00}c_{11}'(\epsilon_1)+A_{2,00}(1-2\epsilon_1) \\
&=&\frac{\epsilon_1-1}{2}\sigma_1^2+\epsilon_1\sigma_0\sigma_2 +\Big(\frac{1}{4}\sigma_1^2+\frac{1}{2}\sigma_0\sigma_2\Big)(1-2\epsilon_1),
\end{eqnarray*}
and therefore the coefficient is 
\[
B_{1,11}=b_{11}'(0)=-\frac{1}{4}\sigma_1^2+\frac{1}{2}\sigma_0\sigma_2 =-B(m).
\]
{\em - Coefficient of $\epsilon_1^2$ in $b_{11}$:} The coefficient is
\[
B_{2,11}=\frac{1}{2}b_{11}''(0)=\frac{1}{2}a_{11}''(0)-A_{2,00}=\frac{1}{4}\sigma_1^2
+\frac{1}{2}\sigma_0\sigma_2-A_{2,00}=0.
\]
{\em - Constant term in $b_{1j}$, $j\geq 2$:} We have
\begin{eqnarray*}
b_{1j}(\epsilon_1)&=&2c_{00}c_{1j}(\epsilon_1)+2c_{01}(\epsilon_1)c_{0j}(\epsilon_1) 
-A_{2,00}(1-\epsilon_j)(1-2\epsilon_1)\\
&=&\epsilon_1(\epsilon_j-1)\sigma_0\sigma_2
+\frac{(\epsilon_1-1)(\epsilon_j-1)}{2}\sigma_1^2 - A_{2,00}(1-\epsilon_j)(1-2\epsilon_1),
\end{eqnarray*}
and therefore the constant term is
\[
B_{0,1j}=b_{1j}(0)=(1-\epsilon_j)\bigg(\frac{\sigma_1^2}{2}-A_{2,00}\bigg)=(1-\epsilon_j)B(m).
\]
{\em - Coefficient of $\epsilon_1$ in $b_{1j}$, $j\geq 2$:} The coefficient is
\[
B_{1,1j}=b_{1j}'(0)=(\epsilon_j-1)\sigma_0\sigma_2+ \frac{\epsilon_j-1}{2}\sigma_1^2+2A_{2,00}(1-\epsilon_j)=0.
\]
We obsrerve that $B_{0,1j}=-(1-\epsilon_j)B_{1,11}$, $j\geq 2$. Hence part \ib of Lemma \ref{lem:vil} gives
\begin{equation}
\epsilon_1 A_{1,11}\Gamma_{11}+\sum_{j=2}^rA_{0,1j}\Gamma_{1j} =O(1),
\label{p3}
\end{equation}
uniformly in $\epsilon_2,\ldots,\epsilon_r$.
Combining (\ref{alv}) and  (\ref{p3}) we conclude that after letting $\epsilon_1\searrow 0$ we are left with
\begin{equation}
I_{r-1}[u]=\sum_{2\leq i\leq j\leq r}b_{ij}\Gamma_{ij} +O(1)\; , \qquad\quad (\epsilon_0=\epsilon_1=0),
\label{alv1}
\end{equation}
uniformly in $\epsilon_2,\ldots,\epsilon_r$. Note that we have the
same coefficients $b_{ij}$ as in (\ref{alv}), unlike the case where the limit $\epsilon_0\searrow 0$ was taken,
in which case we passed from the original coefficients $a_{ij}$ to the coefficients $b_{ij}$.

We proceed in this way. At the $i$th step we denote by $B_{k,ij}$ the coefficient of $\epsilon_i^k$ in $b_{ij}$, $j\geq i$, and observe
that (exactly as in the case $i=1$) there holds
\[
B_{0,ij}=-(1-\epsilon_j)B_{1,ii}\;\quad  , \qquad\quad B_{2,ii}=B_{1,ij}=0\; , \qquad j\geq i+1.
\]
Hence $\ib$ of Lemma \ref{lem:vil} implies the cancelation (modulo uniformly bounded terms)
of all terms that, individually, diverge as $\epsilon_i\searrow 0$. Eventually, after letting $\epsilon_{r-1}\searrow 0$,
we arrive at
\begin{equation}
I_{r-1}[u]=b_{rr}\Gamma_{rr}+O(1) \; , \quad (\epsilon_0=\epsilon_1=\ldots =\epsilon_{r-1}=0),
\label{alvr}
\end{equation}
where $b_{rr}$ has been defined in (\ref{alv}). We observe now that
\[
\int_0^1\frac{u^2}{t^2}X_1^2\ldots X_r^2dt=\Gamma_{rr}.
\]
Hence, using the fact that $\Gamma_{rr}\to +\infty$ as
$\epsilon_r\searrow 0$ (cf (\ref{diverge})) we obtain
\begin{eqnarray*}
\inf_{C^{\infty}_c(0,1)}\frac{I_{r-1}[v]}{\int_0^1\frac{v^2}{t^2}X_1^2\ldots X_r^2dt}
&\leq&\lim_{\epsilon_r\to 0+}\frac{b_{rr}\Gamma_{rr}+O(1)}{\Gamma_{rr}}\\
&=& \lim_{\epsilon_r\to 0+}a_{rr}\\
&=&\lim_{\epsilon_r\to 0} ( c_{0r}^2+2c_{00}c_{rr}) \\
&=&\frac{1}{4}\sigma_1^2-\frac{1}{2}\sigma_0\sigma_2 \\
&=& B(m).
\end{eqnarray*}
This proves part $\ib$ of the theorem. Part $\ia$ follows from (\ref{alvr}) by slightly varrying
the above argument. $\hfill //$

%&&&&&&&&&&&&&&&&&&&&&&&&&&&&&&&&&&&&&&&&&&&&&&&&&&&&&&&&

%\newpage

G. Barbatis \\  
Department of Mathematics \\
University of Ioannina\\ 
45110 Ioannina \\
Greece

\end{document}